\documentclass[11pt]{article}
\usepackage{amssymb}
\usepackage{amsfonts}
\usepackage[latin1]{inputenc}

\newtheorem{theorem}{Theorem}

\newtheorem{definition}[theorem]{Definition}

\newtheorem{proposition}[theorem]{Proposition}

\newenvironment{proof}[1][Proof]{\noindent\textbf{#1.} }{\ \rule{0.5em}{0.5em}}

\begin{document}

\title{KMS States, Entropy and a Variational Principle for Pressure}
\author{Gilles G. de Castro * \, and\, Artur O. Lopes \thanks{%
Instituto de Matem\'atica, UFRGS, 91509-900 Porto Alegre, Brasil. Partially
supported by CNPq, PRONEX -- Sistemas Din\^amicos, Instituto do Mil\^enio,
and beneficiary of CAPES financial support.}}
\date{22 december 2008}
\maketitle

\begin{abstract}
We want to relate the concepts of entropy and pressure to that of KMS states for $C^*$-Algebras.
Several different definitions of entropy are known in our days. The one
we describe here is quite natural and extends the usual one for Dynamical
Systems in Thermodynamic Formalism Theory. It has the advantage of been very
easy to be introduced. It is basically obtained from transfer operators (also called Ruelle operators).
Later we introduce the concept of pressure as a min-max principle.

Finally, we consider the concept of a KMS state as an equilibrium state
for a potential (in the context of $C^*$-Algebras) and we show that there is a relation between KMS states for certain algebras coming from a continuous transformation and equilibrium measures.
\end{abstract}

\newpage

\section{\protect\bigskip Introduction and Main Result}

We want to relate equilibrium measures in the theory of Thermodynamic Formalism with KMS states (the analog of equilibrium measures) in the $C^*$ Algebras theory.
Several different definitions of entropy are known in our days \cite{Ne}.
The one presented here is quite natural and extends the usual one for
Dynamical Systems in Thermodynamic Formalism Theory \cite{L1}. It has the
advantage of been very easy to be introduced. It is basically obtained from
transfer operators.

Later, in section 3, we interpert the concepts of entropy and pressure in the setting of commutative $C^*$-Algebra.

Finally, we consider the concept of a KMS state as an equilibrium state
for a potential (in the context of $C^*$-Algebras) and show that they are related to the equilibrium measures of the Thermodynamic Formalism Theory.

In the next section we will describe in a brief way the main pre-requisites for the
statement of our main result.

\section{\protect\bigskip A brief review of Thermodynamic Formalism and $C^*$-Algebras}

We will present first the main concepts of the theory of Thermodynamic Formalism.
D. Ruelle and Y. Sinai are the more important names in the foundations of this mathematical theory which was inspired in problems borrowed from Statistical Mechanics.

We denote $C(X)$ the space of continuous functions on $X$
taking values on the real numbers where $(X,d)$ is a compact metric space.

Consider the Borel sigma-algebra ${\cal B}$ over $X$ and
a continuous transformation $T:X \to X$. Denote by
${ \cal M} (T)$ the set of invariant probabilities $\nu$ for $T$ (that is, $\int f \circ T \, d \nu= \int f \, d \nu$, for any $f\in C(X)$).
We assume that $T$ is an expanding map (see definition in \cite{R2}).

We refer the reader to \cite{R1} \cite{R2} \cite{R3} for general definitions and properties
of Thermodynamic Formalism and expanding maps.

Typical examples of such transformations (for which that are a lot
of nice results [R2]) are the shift transformation $T$ in the Bernoulli space $\Omega=\{1,2,...,d \}^\mathbb{N}$ and
also  $C^{1+\alpha}$-transformations of the circle such that $|T
'(x)|> c > 1$, where $| \, \, \, |$ is the usual norm (one can
associate the circle to the interval $[0,1)$ in a standard way)
and $c$ is a constant.

The geodesic flow in compact constant negative curvature surfaces
induces in the boundary of Poincar\'e disk a Markov transformation
$G$ such that for some $n$, we have $G^n=T$, and where $T$ is
continuous expanding and acts on the circle (see \cite{BS} \cite{LT3}). Our
results can be applied for such $T$.

We denote by ${\cal H}= {\cal H}_\alpha$ the set of
$\alpha$-Holder functions taking complex values, where
$\alpha$ is fixed $0< \alpha\leq 1$.

For each $\nu\in { \cal M} (T)$, $h(\nu)$ denotes the
the Shannon-Kolmogorov entropy of $\nu$ (see \cite{PP} for definition).

The entropy measures the dynamic complexity of the action of the transformation $T$ in sets of measure one.

One problem one can be interested is maximizing entropy among invariant probabilities.

That is consider,
$h(T)= \sup \{h(\nu)\,|\,  \nu \in { \cal M} (T)\}$. The value
$h(T)$ is called the topological entropy
of $T$. A probability that attains such supremum value is called a measure of maximal entropy.

In this way we are looking for the probability with the largest complexity. In the case of the Bernoulli space
$\{1,2,...,d \}^\mathbb{N}$ the maximal value of the entropy of invariant probabilities is $\log (d)$ and there is a unique probability $\mu$
which attains such value. In this case the maximal entropy measure $\mu$ is the independent probability with weights $1/d$.

In physics one of the main principles is that Nature maximizes entropy. The above maximal entropy measure   $\mu$  corresponds in Statistical Mechanics to what is expected  at temperature infinite (see \cite{R1} \cite{R2}).   When the temperature is finite Nature maximizes Pressure.  A probability that maximizes pressure is called a Gibbs probability (or state). There exists in fact an external potential $A: \Omega \to \mathbb{R}$ which describe the interaction $A(w) = A(w_0,w_1,w_2,...) $, where $ w =(w_0,w_1,w_2,...)\in \Omega$,   around neighborhoods in the lattice $\Omega= \{1,2,...,d \}^\mathbb{N}$. The simplest case is when this potential is a function $A(w) = A(w_0,w_1,w_2,...)$ which depends only in a finite number of coordinates. For example, if depends on two coordinates
then $A(w) = A(w_0,w_1)$. In this case we have a finite range iteration potential. These  kind of potentials are more easy to deal, but the more important potentials, for mathematical applications, are the ones which depends on the
whole $ w=(w_0,w_1,w_2,...) \in  \{1,2,...,d \}^{\mathbb{N}}$.

A common example in Statistical Mechanics is when $d=2$ and one considers that $1$ corresponds to the spin $+$ and $2$ corresponds to the spin $-$. An element $w$ in Bernoulii space could be $w=(+\,-\,-\,+\,-\,+\,-\,+\,+\,...)$, which means an element $w$  where the spin is up or down in different positions in a lattice over the set $\mathbb{N}$. The Gibbs probability for $A$
describes probabilities of Borel sets of in the space $\{+,-\}^{\mathbb{N}}$ which are determined by the interactions that are given by $A$.

We can also consider, as in Statistical Mechanics, an extra real parameter $\beta$ which represents $1/T$, where $T$ is temperature.

\begin{definition} Given a potential $A$, the pressure of $A$ at temperature $T=1/\beta$ is, by definition,
$$ P(\beta\, A)  =\sup \{h(\nu)\, +\, \beta \,\int A d \nu\,|\,  \nu \in { \cal M} (T)\}.$$
\end{definition}

A measure $\mu=\mu_A$ satisfying such supremum is called a Gibbs state for $A$ at temperature $T$. It describes what is physically observed in probabilistic terms (see \cite{R1}). The probability $\mu_A$ is also called an equilibrium state for $A$.

When $\beta=1$ we just say $P(A)$.

If $A$ is Holder the Gibbs state $\mu_A$ is unique (see \cite{R3} \cite{PP}). The Holder property of the potential $A$ corresponds to a rapidly decay of interaction between neighborhoods.

In the case $\beta=0$ (which means $A=0$) we get the case we described before (temperature infinite). In this case the Gibbs  state is the independent probability we previously consider.

For a differentiable transformation $T: S^1 \to S^1$, where $S^1$ is the circle (or, an interval), a  very important potential to be  consider is
$A(x)\,=\, - \log T' (x)$. In this case $-\int \log T' (x) d \mu(x)$ measures the $\mu$-mean sensibility with respect to initial conditions. In this case
the measure $\mu$ which maximizes pressure is called Bowen-Ruelle-Sinai probability (see \cite{MS}).

For more general transformations on the circle one can consider the extra parameter $\beta$ (which now has nothing to do with temperature)
and the potentials $\beta (- \log T')$. A special value of $\beta$, namely the one such that $P(- \beta  \log T')=0$, will be associated to the Hausdorff dimension of sets which are important
from the dynamical point of view \cite{R2} \cite{MU}.

Applications of this theory for dimension of fractals, geometry, zeta functions, etc... can be found in \cite{CR} \cite{BS}, \cite{MU}, \cite{ PP}, \cite{LT3}.

The main tool for obtaining the Gibbs probability is the Ruelle operator ${\cal L}_A: C(X) \to C(X)$ (which in Statistical Physics is called the transfer operator).

Let's consider the general case:

\begin{definition}
Given
$A:X\to \mathbb{R}$,  the Ruelle operator ${\cal L}_A: C(X) \to C(X)$ is defined in the following way: given $f:X\to \mathbb{R}$,  we denote $g = {\cal L}_A (f)$, in such way that for a $x\in X$, we have
$g(x)=({\cal L}_A (f))(x) = \sum_{T(z)=x} e^{A(z)}\, f(z)$.

\end{definition}

The Ruelle operator ${\cal L}_A$ is also called here transfer operator.

We can also consider the dual Ruelle operator  ${\cal L}_A^*$ acting on measures over the Borel sigma algebra on $X$.

If $A$ is Holder then the Ruelle operator can also act in the space of Holder functions ${\cal H}_\alpha$.

When $A$ is such that
${\cal L}_A(1) =1$, we say that the potential $A$ is normalized. In this case if $\nu$ is a probability, then ${\cal L}_ A^*(\nu)$ is also a probability.

We will state now a main result in the theory in the particular  case of the Bernoulli space. The next theorem is a more advanced version of Perron Theorem for positive matrices.

\begin{theorem} Ruelle Theorem \cite{R3} \cite{PP}

If $A: \{1,2,...,d \}^\mathbb{N}\to \mathbb{R}$ is Holder, then there exists a maximal eigenvalue $\lambda$  for ${\cal L}_A$ and
a Holder eigenfunction $\phi$ such that

a)  ${\cal L}_A(\phi) = \lambda \, \phi$

Moreover, there exists a eigen-probability $\nu$ such that

b)  ${\cal L}_A^*(\nu) = \lambda \, \nu$.

Finally, the Gibbs state probability $\mu_A$  for $A$ is $\mu_A= \phi \, \nu$ (after suitable normalization).

\end{theorem}

We point out that when  $A$ depends only on two coordinates the above theorem is a consequence of Perron Theorem for positive matrices (see \cite{PP}).

In the case the potential is Holder the main eigenvalue $\lambda$ is isolated in the spectrum of the operator ${\cal L}_A$ (when acting in the space of Holder functions)  \cite{PP}.

In \cite{L1} it is presented a different way to compute entropy via the Perron operator acting on different potentials.

\begin{theorem} Denote $ \mathbb{B}^{+}$ the set of Borel positive functions on $\Omega$.
    Given $\mu \in \mathcal{M}(T)$ and a Holder potential $A$, the  entropy of $\mu$ is given by
   $$h(\mu)=\inf_{f \in \mathbb{B}^{+}} \int \ln\left(\frac{P_{A}f}{A\, f}\right) d\nu.$$
\end{theorem}

This theorem show that one can avoid the dynamical point of view of entropy (of the usual definition which considers partitions of the Bernoulli space, refinements of the partition by iteration of $T$, and so on...) and address all the computation to the action of the Ruelle operator. This result is quite useful for a generalization for the case of $C^*$ Algebras where there is no natural dynamics involved
(or a definition based on dynamics would be quite complicated).

In \cite{L1} it is presented a different way to compute pressure via a min-max principle.

\begin{theorem} Denote $ \mathbb{B}^{+}$ the set of Borel positive functions on $\Omega$.

   For  a Holder potential $A$ the topological pressure is given by the expression
   $$P(A)=
     \sup_{\mu \in \mathcal{M}(T)} \inf_{f \in \mathbb{B}^{+}}
      \int \log(\frac{P_{A}\,f}{f}) d\mu.$$
\end{theorem}
\vspace{0.3cm}

Now we will  briefly describe  some basic results on $C^*$-Algebras. This theory is presented in a quite elegant way   in references \cite{Ped} \cite{Bra}. I.M Gelfand and J. Von Neumann are the more important names in the foundations of this theory.

We refer the reader to  \cite{Ren1} \cite{RE1}
\cite{exel01} \cite{exel02} \cite{exel03} \cite{EL2} for a more complete description of the relation of Thermodynamic Formalism and $C^*$-Algebras.

The KMS states are quite important in Quantum Statistical Mechanics \cite
{Bra}. They play in $C^*$-Algebras  the role of equilibrium states.
We will explain this bellow.

 \vspace{0.3cm}

First of all the space $C(X)$ will denote now the space of continuous functions defined on the compact metric space $X$ taking values on $\mathbb{C}$.

\begin{definition} A complete normed algebra over $\mathbb{C}$ with an involution operation $*$ such that
$$\| a \, a^* \| = \| a\|^2$$
is called a $C^*$-Algebra.
\end{definition}

We refer the reader to \cite{Bra} for more detailed  definitions and main properties of $C^*$-Algebras (see definition 2.1.1 and example 2.1.2).

In Quantum Mechanics the potential (also called observable) $A:X\to \mathbb{C}$ will be replaced by an operator (acting on the complex Hilbert space  ${\cal L}^2(\mu)$). In this way the commutative algebra of functions (with the usual complex product structure) gives place to the non-commutative
algebra of bounded operators $B$ from ${\cal L}^2(\mu)$ to itself (where the product structure is the composition of operators). Moreover, for the $*$ operation on the algebra, we consider for each $B$, the new operator  $B^*$  which means the adjoint of the operator $B$. The norm in the algebra is the operator norm.

For $T$, we consider here two $C^*$-Algebras ${\cal U}$ and ${\cal V}$ which are subalgebras of the algebra of bounded operators $B: {\cal L}^2(\mu)\to {\cal L}^2(\mu)$,
where $\mu$ is a Gibbs measure for a fixed potential $\tilde {A}$. In the case of the Bernoulli space $\Omega=\{1,2,..,d\}^\mathbb{N}$ (the main case we consider here) this potential $\tilde {A}$ can be taken the constant potential $-\log d$ (see \cite{exel03}). In this case, $\mu$ is the independent probability (with weights $1/d$) over $\{1,2,..,d\}^\mathbb{N}$. Moreover,
for the dual of Ruelle operator for $\tilde{A}$  acting on probabilities we have that ${\cal L}_{\tilde{A}}^{*} (\mu)=\mu$ (see \cite{exel03} \cite{PP}).

\begin{definition}
An important class of linear operators are the ones  of the form
$M_f : {\cal L}^2 (\mu) \to {\cal L}^2 (\mu)$, for a given fixed $f\in C(X)$,
and defined by $M_f (\eta) (x)= f(x) \eta (x)$, for any $\eta$ in ${\cal L}^2 (\mu)$.
\end{definition}

In order to simplify the notation, sometimes people denote by $f$ the linear operator $M_f$.

Note that for $M_f $ and $M_g$, $f,g \in C(X)$, $X$ metric compact space, the product operation satisfies
$M_f \circ M_g = M_ {f.g}$, where $\,.\,$ means multiplication
over the complex field $\mathbb{C}$.

Note that the   $*$ operation applied on $M_f$, $f \in C(X)$,
is given by $M_f^{*}= M_{\overline{f}}$, where $\overline{z}$
is the complex conjugated
of  $z\in\mathbb{C}$.
In this way, $M_{\overline{f}}$ is the adjoint operator of $M_f$ over
${\cal L}^2 (\mu)$.

\begin{definition}
Denote by $S: {\cal L}^2 (\mu) \to {\cal L}^2 (\mu)$ the Koopman  operator
where for $ \eta \in {\cal L}^2 (\mu)$, we define
$(S \eta) (x) = \eta (T(x))$. Such $S$ defines a bounded linear operator in
${\cal L}^2 (\mu)$.
\end{definition}

In Thermodynamic Formalism
it is usual to consider the Koopman operator $S$  acting on ${\cal L}^2 (\mu)$, and it is well known that
its adjoint (over ${\cal L}^2 (\mu)$) is the operator
${\cal L}_{\tilde{A}}$ acting on ${\cal L}^2 (\mu)$ (it is well defined as one can see in \cite{PP}).

The main point for our choice of $\mu$ is exactly that ${\cal L}_{\tilde{A}}=S^*$.

Now we have all elements to define our $C^*$-Algebras.

\begin{definition} For $T$ and a Gibbs measure $\mu$, we define $\cal V=\cal V(T,\mu)$ as the sub-$C^*$-Algebra of the $C^*$-Algebra of bounded linear operators  on ${\cal L}^2 (\mu)$ generated by $S$ and $M_f$ for all $f\in C(X)$. We also define the $C^{*}$-Algebra ${\cal U}= {\cal U} (\mu,T)$ as being generated by the elements of the form
$M_f S^n (S^{*})^n M_g$, where $n\in \mathbb{N}$ and $f,g\in C(X)$.

The algebra $\cal U$ is a sub-$C^*$-Algebra of $\cal V$. We can see that each element $a$ in ${\cal V}$ is the limit of finite
sums  $\sum_{i} M_{f_i} S^{n_i} (S^{*})^{m_i} M_{g_i}$, whereas an element $a$ in ${\cal U}$ is the limit of finite
sums  $\sum_{i} M_{f_i} S^{n_i} (S^{*})^{n_i} M_{g_i}$ (the exponents of $S$ and $S^{*}$ are the same).

\end{definition}

\begin{definition}

An element $a$ in a $C^{*}$-Algebra is positive, if it is of the
form $a=b b^{*}$ with $b$ in the $C^{*}$-Algebra.

\end{definition}

\begin{definition}
A state in a $C^*$-Algebra with unit $A$ is a linear functional $\phi: A \to \mathbb{C}$
such that

a) $\phi (1)=1$

b) $\phi (a) $ is a positive real number
for each positive element      $a$ on the $C^{*}$-Algebra
$A$.

\end{definition}

A state $\phi$ in the context of $C^{*}$-Algebras plays the role of a probability $\nu$ in Thermodynamic Formalism.

\begin{definition} An one-parameter group of automorphisms in a $C^{*}$-Algebra $A$ is a strongly continuous group homomorphism $\sigma:\mathbb{R}\to \mathrm{Aut}(A)$. We denote the automorphism $\sigma(t)$ by $\sigma_t$.
\end{definition}

We think $\sigma$ as a dynamic temporal
evolution in the $C^{*}$-Algebra $A$.

\begin{definition}
An element $a\in A$ is called analytic for $\sigma$, if
$\sigma_t(a)$ has an analytic extension from $t\in\mathbb{R}$ to all
$z\in \mathbb{C}$.

\end{definition}

The set of analytical elements is always dense in $A$.

\begin{definition} Given $\beta\in \mathbb{R}$ and an one-parameter group of automorphisms, a state $\phi$ in $A$ is a $(\sigma,\beta)$-KMS state if for any $b\in A$ and any analytical element $a\in A$, we have

$$ \phi (a . b) = \phi (b . \sigma_{\beta i} (a)).$$

\end{definition}

We now consider certain dynamical evolutions in the $C^{*}$-Algebras $\cal U$ and $\cal V$.

\begin{definition}

Given a strictly positive function $H:X\to \mathbb{R}  $ we define an associated one-parameter group of automorphisms
$\sigma:\mathbb{R}\to\mathrm{Aut}(\cal V)$, where for each  $t\in \mathbb{R}$ we have that
$\sigma_t$ is defined by:

a) for each fixed $t\in  \mathbb{R}$ and any $M_f$, we have
$\sigma_t (M_f)=M_f$,

b) for each fixed $t\in  \mathbb{R}$, we have
$\sigma_t (S)= M_{H^{i  t}} \circ S$,
, in the sense that
$(\sigma_t  (S) (\eta)) (x) = H^{i t}(x) \eta(T(x))\in
{\cal L}^2 (\mu),$ for any $\eta\in
{\cal L}^2 (\mu)$.

Since for each  $t\in \mathbb{R}$ we have that $\sigma_t(\cal U)\subset\cal U$, we can restrict $\sigma$ to $\cal U$.
\end{definition}

In terms of the formalism of  $C^{*}$-dynamical systems,
the positive function $H$ defines the dynamics
of the evolution with time $t\in \mathbb{R}$ of a  $C^{*}$-dynamical system.
The dynamics of the shift transformation $T$ plays here just the role of spatial translation on the lattice.

The $H$ above corresponds to  the $A$ we consider before in Thermodynamic Formalism by the expression $H=e^A$. If we introduce a parameter
$\beta$, then we have to consider in the $C^*$-Algebra setting the potential $H^\beta$.
Our purpose here is to analyze such system  for each pair $(H,\beta)$.

For $H$ and $\beta$ fixed (that is, we consider the potential $H^\beta$),
we denote a KMS state by $\phi_{H,\beta}$
and we leave $\phi$ for a general $C^{*}$-dynamical system state.

The state $\phi_{H,\beta}$ is what is expected from the Quantum Statistical point of view of a system  governed by $H$ under temperature $T=1/\beta$ (see\cite{Bra}).

It is easy to see  that for $H$ and $\beta$ fixed, the condition
$$ \phi (a . b) = \phi (b . \sigma_{\beta i} (a)),$$
is equivalent to $\forall \tau\in \mathbb{R}$,
$$ \phi (\sigma_{\tau}(a) . b) = \phi (b . \sigma_{\tau + \beta i} (a)).$$

It follows from section 8.12 in \cite{PP} that if $\phi$ is a KMS state
for $H,\beta$, then for any analytic $a\in {\cal U}$, we
have that the extension
$\tau\to \phi(\sigma_\tau ( a) ) $ to $z\to \phi(\sigma_z ( a) ) $ is
a bounded entire function
and therefore constant. In this sense $\phi$ is stationary.

A natural question is: for a given $\beta$ and $H$,
when the KMS state $\phi_{H,\beta}$ exist and when it  is unique?
This question is consider in \cite{exel02} for the $C^*$-Algebra $\mathcal{V}$ and in \cite{exel03} for the  $C^*$-Algebra $\cal U$.
Another presentation of the uniqueness part of this result appears in \cite{EL2}.

Note that the action of the linear functional $\phi$ in the set of operators $M_f$, where $f$ range in all continuous functions,
defines, via Riesz Theorem, a measure $\nu$ over $X$. That is, for all $f \in C(X)$ we have $\phi(M_f)= \int f \, d \nu$. In fact $\nu$ is a probability from the hypothesis we assume for a $C^*$-state $\phi$.

One of the main points in \cite{exel03} is that for the KMS state in $\cal U$ associated to $H$ and $\beta$ this measure $\nu$ is the eigen-measure $\nu_{H,\beta}$ for the dual Ruelle operator ${\cal L}_{A}^*$, where $A=-\beta\, \log H$, we consider before.

In this way, one can associate, in a unique way,
each KMS state $\phi_{H,\beta}$ to the eigen-measure $\nu_{H,\beta}$. This is an interesting relation of Thermodynamic Formalism and $C^*$-Algebras.

For the case of KMS states in $\cal V$, we have the extra condition that the pressure of $H^{-\beta}$ is zero so that we have the existance of KMS states for only on value of $\beta$.

\section{\protect\bigskip Proof of the main result}

Suppose $A$ is a commutative C*-Algebra with unity and $\alpha :A\rightarrow
A$ is an injective endomorphisms which preserves unity. Considering $A=C(X)$%
, via Gelfand-Naimark theorem, we get that $\alpha $ is of the form $\alpha
(a)=a\circ T$, where $T:X\rightarrow X$ is a continuous transformation.

\begin{definition}
A transfer operator for $\alpha $ is linear transformation $%
L:A\rightarrow A$, such that $L(\alpha (a)b)=aL(b)$, $\forall a,b\in A$.
Moreover, if $L(1)=1$, then we say that the transfer operator is normalized.
\end{definition}

In the case $A$ is a commutative algebra, the transfer operator takes the
form of a Ruelle operator \cite{PP} \cite{DutJorgDisint} \cite{DutkayJorg}
\cite{Ren1} \cite{exel01}.

\begin{proposition}
\cite{kwas01} If $L$ is a transfer operator for $\alpha $ then:

\begin{enumerate}
\item $L(a\alpha (b))=L(a)b$, $\forall a,b\in A$;

\item $L(1)$ is a central positive element of de $A$.
\end{enumerate}
\end{proposition}

\begin{proposition}
\cite{kwas01} If $T$ is a local homeomorphism, then all transfer
operators for $\alpha $ are of the form:%
\[
L(a)|_{x}=\sum_{T(y)=x}\rho (y)a(y),
\]%
where $\rho :X\rightarrow \lbrack 0,\infty )$ is a continuous function.
Moreover, any continuous function $\rho :X\rightarrow \lbrack 0,\infty )$,
defines a transfer operator by the above expression.
\end{proposition}

 In the notation of last section $\rho=e^A$.

From now on, we suppose $T$ is a local homeomorphism and we denote by $%
L_{\rho }$ transfer operator defined by $\rho $. We want to introduce a
notion of entropy for a state $\phi :A\rightarrow \mathbb{C}$, using the
transfer operator. This generalizes the point of view described in \cite{L1}
\cite{L2}. Denote $A_{+}:=\{a\in A~|~\sigma (a)\in (0,\infty )\}$

\begin{definition}
Given a state $\phi $ in $A$, we define the entropy of $\phi $ by:%
\[
h(\phi )=\inf_{a\in A_{+}}\left\{ \phi \left( \ln \left( \frac{L_{\rho }(a)}{%
\rho a}\right) \right) \right\},
\]%
where $L_{\rho }$ is the transfer operator for $\rho :X\rightarrow (0,\infty
)$.
\end{definition}

Note that the above definition is independent of the choice of $\rho $,
indeed, if $\rho ^{\prime }:X\rightarrow (0,\infty )$ is another continuous
function, then taking $a^{\prime }=a\rho (\rho ^{\prime })^{-1}$, for some
arbitrary $a\in A$, we get that $a^{\prime }\in A_{+}$, and%
\[
\frac{L_{\rho ^{\prime }}(a^{\prime })}{\rho ^{\prime }a^{\prime }}=\frac{1}{%
\rho a}\sum_{y=T(x)}\rho ^{\prime }(x)a(x)\rho (x)\rho ^{\prime }(x)^{-1}=%
\frac{L_{\rho }(a)}{\rho a}.
\]

Therefore we are considering the infimum over the same set.

\begin{definition}
We say that a state $\phi $ is $\alpha $-invariant, if, $\phi \circ \alpha
=\phi $.
\end{definition}

\begin{definition}
Given an element $b\in A_{+}$, we define the topological pressure of $b$ by%
\[
p(b)=\sup_{\phi ~\mathrm{inv}}\left\{ h(\phi )+\phi (\ln b)\right\} .
\]%
If $\phi $ is an $\alpha$-invariant state, such that, $h(\phi )+\phi (\ln
b)=p(b)$, Then, we say that $\phi $ is an $C^*$-equil\'{\i}brium state for $%
b $.
\end{definition}

\begin{proposition}
If $L_{\rho }(1)=1$, then, there exists an state $\phi $ such that $\phi
\circ L_{\rho }=\phi $.
\end{proposition}

\begin{proof}
As $L_{\rho }(1)=1$, we have that $L_{\rho }^{\ast }(\mathcal{S)\subset S}$,
where $\mathcal{S}$ is the set of all states of $A$. Using
Tychonoff-Schauder theorem, we get that $L_{\rho }^{\ast }|_{\mathcal{S}}$
has a fixed point.
\end{proof}

\begin{proposition}
\label{TranfNorm} If $L_{\rho }(1)=1$, then $p(\rho )=0$. Moreover, the
states $\phi $ which satisfies $\phi \circ L_{\rho }=\phi $ are equilibrium
states for $\rho $.
\end{proposition}

\begin{proof}
Note that using $L_{\rho }$ in the definition of entropy we have that%
\[
p(\rho )=\sup_{\phi ~\mathrm{inv}}\left\{ \inf_{a\in A_{+}}\left\{ \phi
\left( \ln \left( \frac{L_{\rho }(a)}{a}\right) \right) \right\} \right\} .
\]

Choosing $a=1$, inside the infimum, it follows that $p(\rho )\leq 0$.

By the other hand, as $L_{\rho }$ is normalized, we have that $L_{\rho
}\circ \alpha =Id$, and, if $\phi \circ L_{\rho }=\phi $, then $\phi \circ
\alpha =$ $\phi \circ L_{\rho }\circ \alpha =\phi $. As $\ln $ is a concave
funtion, then $\ln (L_{\rho }(a))\geq L_{\rho }(\ln a)$. Therefore,
\[
\phi \left( \ln \left( \frac{L_{\rho }(a)}{a}\right) \right) =\phi \left(
\ln (L_{\rho }(a))-\ln a\right) \geq \phi \left( L_{\rho }(\ln a)-\ln
a\right) .
\]

If $\phi \circ L_{\rho }=\phi $, the right hand side of the inequality above
is equal to zero, and therefore, $\inf_{a\in A_{+}}\left\{ \phi \left( \ln
\left( \frac{L_{\rho }(a)}{a}\right) \right) \right\} =0$. It follows that $%
p(\rho )\geq 0$, and in the case of an eigen-state, we have $h(\phi )+\phi
(\rho )=\inf_{a\in A_{+}}\left\{ \phi \left( \ln \left( \frac{L_{\rho }(a)}{a%
}\right) \right) \right\} =0=p(\rho )$.
\end{proof}

\bigskip

If we consider the context of an algebra $A$, an injective endomorphism
preserving unity $\alpha $, and a normalized transfer operator $L$, we can
consider (among others) two different C*-algebras: the cross-product
endomorphism $A\rtimes _{\alpha ,L}\mathbb{N}$ (see \cite{exel01}) and the
C*-algebra given by approximately proper equivalence relations $C^{\ast
}(\mathcal{R},\mathcal{E})$ (see \cite{exel03}). The second algebra is
related with the equivalence relation $x\sim y\iff \exists n\in \mathbb{N}$,
such that, $T^{n}(x)=T^{n}(y)$. The first one considers a more broad
equivalence relation $x\sim y\iff \exists n,m\in \mathbb{N}$, such that, $%
T^{n}(x)=T^{m}(y) $.

The algebra $C^{\ast}(\mathcal{R},\mathcal{E})$ is a generalization of the algebra $\cal U$ in the previous section whereas $A\rtimes _{\alpha ,L}\mathbb{N}$ is a generalization of the algebra $\cal V$. In fact in the context of the previous section, for each Gibbs measure $\mu$ we find representations of $A\rtimes _{\alpha ,L}\mathbb{N}$ and $%
C^{\ast }(\mathcal{R},\mathcal{E})$ in the Hilbert space $\mathcal{L}^2(\mu)$ such that the images are $\mathcal{V}$ and $\mathcal{U}$ respectively.

We want to relate KMS states of $A\rtimes _{\alpha ,L}\mathbb{N}$ and $%
C^{\ast }(\mathcal{R},\mathcal{E})$, with the equilibrium states (in $A$) of
the potential $h^{-\beta }$. Here $\beta $ represents the inverse
of temperature.

In the case the algebra $A\ $ is commutative, we have unique conditional
expectations $F:A\rtimes _{\alpha ,L}\mathbb{N\rightarrow }A $ and $%
G:C^{\ast }(\mathcal{R},\mathcal{E})\rightarrow A$. Moreover, if $E:=\alpha
\circ L:A\rightarrow \alpha (A)$, for a conditional expectation with finite
index, then the KMS states $\psi $ of $A\rtimes _{\alpha ,L}\mathbb{N}$ can
be decomposed as $\psi =\phi \circ F$, where $\phi $ is a state on $A$ which
satisfies%
\[
\phi (a)=\phi (L(\Lambda a)),~\forall a\in A,
\]%
and $\Lambda =h^{-\beta }\mathrm{ind}(E)$. The KMS state $\psi $ of $C^{\ast
}(\mathcal{R},\mathcal{E})$ can be decomposed as $\psi =\phi \circ G$, where
$\phi $ is a state in $A$ which satisfies%
\[
\phi (a)=\phi (\Lambda ^{-[n]}E_{n}(\Lambda ^{\lbrack n]}a)),~\forall a\in
A~,\forall n\in \mathbb{N},
\]%
where $E_{n}=\alpha ^{n}\circ L^{n}$, and $\Lambda ^{\lbrack
n]}=\prod_{i=0}^{n-1}\alpha ^{i}(h^{-\beta }\mathrm{ind}(E))$.

\begin{proposition}
If $\psi =\phi \circ F$ is a $(h,\beta )$-KMS state for $A\rtimes _{\alpha
,L}\mathbb{N}$, and $L(\Lambda 1)=1$, then $\phi $ is an equilibrium state
(in $A$) for the potential $h^{-\beta }$.
\end{proposition}

\begin{proof}
The condition $L(\Lambda 1)=1$ implies that $L_{h^{-\beta }}$ is a
normalized transfer operator, and, therefore $p(h^{-\beta })=0$. Moreover,
the KMS condition says that $\phi (a)=\phi (L_{h^{-\beta }}(a))$, which
implies that $\phi (\alpha (a))=\phi (a)$. It follows from Proposition \ref%
{TranfNorm} that $\phi $ is an equilibrium state for $h^{-\beta }$.
\end{proof}

\bigskip

In the construction of the algebras we are interested, the choice of two
normalized transfer operator define isomorphic algebras, in such way that,
the choice of the operator can be arbitrary.

Suppose that $L_{\rho }(k)=\lambda k$, for some $\lambda >0$ and $k\in A_{+}$%
. Defining $\widetilde{\rho }=\frac{\rho k}{\lambda \alpha (k)}$, we have
that $L_{\widetilde{\rho }}$ is a normalized transfer operator for $\alpha $%
, and, therefore, we can use it to obtain $C^{\ast }(\mathcal{R},\mathcal{E}%
) $.

\begin{proposition}
Suppose that $\psi =\phi \circ G$ is a $(h,\beta )$-KMS state of $C^{\ast }(%
\mathcal{R},\mathcal{E})$. Consider $\rho =h^{-\beta }$ and suppose that $%
L_{\rho }(k)=\lambda k$, for some $\lambda >0$ e $k\in A_{+}$. Denote $%
\widetilde{\rho }=\frac{\rho k}{\lambda \alpha (k)}$, and consider $%
\widetilde{\phi }$ the state of $A\ $ given by $\widetilde{\phi }(a)=\phi
(ka)$. If $\lim_{n\rightarrow \infty }\left\Vert L_{\widetilde{\rho }%
}^{n}(a)-\widetilde{\phi }(a)\right\Vert =0$, $\forall a\in A$, then $%
\widetilde{\phi }$ is an equilibrium state for $\widetilde{\rho }$.
\end{proposition}

\begin{proof}
We can suppose, w.l.g., that $C^{\ast }(\mathcal{R},\mathcal{E})$ was
obtained from $L_{\widetilde{\rho }}$, in such way that $\mathrm{ind}(E)=%
\widetilde{\rho }^{-1}$. Therefore,
\[
\Lambda ^{\lbrack 1]}=\left( \rho \widetilde{\rho }^{-1}\right) =\rho \frac{%
\lambda \alpha (k)}{\rho k}=\frac{\lambda \alpha (k)}{k},
\]%
and, more generally%
\[
\Lambda ^{\lbrack n]}=\prod_{i=0}^{n-1}\alpha ^{i}\left( \rho \widetilde{%
\rho }^{-1}\right) =\lambda ^{n}\frac{\prod_{i=0}^{n-1}\alpha ^{i+1}(k)}{%
\prod_{i=0}^{n-1}\alpha ^{i}(k)}=\frac{\lambda ^{n}\alpha ^{n}(k)}{k}
\]%
The KMS condition implies%
\[
\phi (a)=\phi \left( \frac{k}{\lambda ^{n}\alpha ^{n}(k)}\alpha ^{n}L_{%
\widetilde{\rho }}^{n}\left( \frac{\lambda ^{n}\alpha ^{n}(k)}{k}a\right)
\right) =
\]%
\[
=\phi \left( k\alpha ^{n}L_{\widetilde{\rho }}^{n}\left( \frac{a}{k}\right)
\right) ,
\]%
for all $n\in \mathbb{N}$. It follows that%
\[
\widetilde{\phi }(a)=\phi (ak)=\phi \left( k\alpha ^{n}L_{\widetilde{\rho }%
}^{n}\left( \frac{a}{k}k\right) \right) =\widetilde{\phi }\left( \alpha
^{n}L_{\widetilde{\rho }}^{n}\left( a\right) \right) .
\]%
Now,%
\[
\left\vert \widetilde{\phi }(L_{\widetilde{\rho }}(a)-a)\right\vert
=\left\vert \widetilde{\phi }(\alpha ^{n}(L_{\widetilde{\rho }}^{n+1}(a)-L_{%
\widetilde{\rho }}^{n}(a))\right\vert \leq
\]%
\[
\leq \widetilde{\phi }\left( \left\Vert \alpha ^{n}(L_{\widetilde{\rho }%
}^{n+1}(a)-L_{\widetilde{\rho }}^{n}(a))\right\Vert \right) \leq \widetilde{%
\phi }\left( \left\Vert L_{\widetilde{\rho }}^{n+1}(a)-L_{\widetilde{\rho }%
}^{n}(a)\right\Vert \right)
\]%
\[
\leq \widetilde{\phi }\left( \left\Vert L_{\widetilde{\rho }}^{n+1}(a)-%
\widetilde{\phi }(a)\right\Vert \right) -\widetilde{\phi }\left( \left\Vert
\widetilde{\phi }(a)-L_{\widetilde{\rho }}^{n}(a)\right\Vert \right)
\longrightarrow _{n\rightarrow \infty }\,\,0,
\]%
and, therefore $\widetilde{\phi }\circ L_{\widetilde{\rho }}=\widetilde{\phi
}$. From proposition \ref{TranfNorm}, the claim follows.
\end{proof}

\bigskip

Note that the hypothesis of the convergence of $L_{\widetilde{\rho }}^{n}$ is
one of the conclusions of Ruelle-Perron-Frobenius theorem (see \cite{PP} \cite{exel03}, \cite{DutkayJorg})
 so that the
classical setting satisfies the hypothesis of the previous proposition.

\end{document}